\documentclass[a4paper, 12pt]{amsart}

\usepackage{amssymb}

\mathchardef\emptyset="001F

\theoremstyle{plain}

\newtheorem{theorem}{Theorem}[section]

\newtheorem{remark}[theorem]{Remark}

\theoremstyle{definition}
\theoremstyle{remark}
\numberwithin{equation}{section}

\newcommand{\Ao}{{\mathbb A}}

\newcommand{\NN}{{\mathbb N}}

\newcommand{\T}{{\mathbb T}}

\newcommand{\V}{{\mathcal V}}

\renewcommand{\div}{\hbox{{\rm div}}}
\newcommand{\grad}{\hbox{{\rm grad}}}

\def\2s{\stackrel{\rm 2s}{\rightharpoonup}}

\title[]
{A bound on group velocity for Bloch wave packets}
\author[Gr\'egoire Allaire]{Gr\'egoire Allaire$^{1}$}
\author[Mariapia Palombaro]{Mariapia Palombaro$^{2}$}
\author[Jeffrey Rauch]{Jeffrey Rauch$^{3}$}

\begin{document}

\maketitle

{
\small
\noindent
$^1$ Centre de Math\'ematiques Appliqu\'ees, \'Ecole Polytechnique, 
91128 Palaiseau, France.\\
Email: gregoire.allaire@polytechnique.fr\\
\noindent
$^2$ 
Institute for Mathematics, University of W\"urzburg,
Am Hubland,
97074 W\"urzburg, Germany.\\
Email: mariapia.palombaro@mathematik.uni-wuerzburg.de \\
\noindent
$^3$ Department of Mathematics, University of Michigan,
Ann Arbor 48109 MI, USA.\\
Email: rauch@umich.edu
}

\maketitle

\section{Main result}

This short note is a sequel to our previous papers \cite{apr}, \cite{apr2} 
on the asymptotic behavior of Bloch wave packet solutions of the wave equation 
in periodic media. 
The purpose is to prove that the group velocity for these Bloch wave packets 
is bounded by the maximal speed of propagation for the original wave equation. 
This follows from the fact 
that the wave packets
provide accurate approximate solutions.
What was lacking is a  mathematical proof that 
uses only the 
definition of the group velocity. 
We follow the notations in \cite{apr2}. 

For periodic coefficients $A_0(y)$ and $\rho_0(y)$ in $L^\infty(\T^N)$, consider the Bloch spectral cell problem
\begin{equation}
\label{celleq}
- (\div_y +2i\pi\theta) \Big( A_0(y) (\grad_y +2i\pi\theta)\psi_n \Big) 
\ =\
 \lambda_n(\theta) \rho_{0}(y) \, \psi_n  \quad \mbox{ in }\ \  \T^N \,,
\end{equation}
with the Bloch parameter $\theta\in [0,1[^N$ and $\T^N$ the unit torus. 
Assuming that $A_0(y)$ is symmetric and uniformly coercive and that $\rho_0(y)$ 
is uniformly bounded away from zero, it is well known that \eqref{celleq} 
admits a countable infinite family of positive real eigenvalues $\lambda_n(\theta)$ 
(repeated according to their multiplicity) 
and associated eigenfunctions $\psi_n(\theta,y)$ which, as functions 
of $y$, belong to $H^1(\T^N)$ \cite{bloch,blp,reedsimon,wilcox}. 
The eigenvalues, being labeled by increasing order, are Lipschitz functions 
of $\theta$ (not more regular because of possible crossings). 
However, simple eigenvalues are analytic functions of $\theta$ \cite{kato}. 
Being a simple eigenvalue is a generic property \cite{albert}.
Normalize the eigenfunctions by
\begin{equation}
\label{eq:psinorm}
\int \rho_0(y)\, |\psi_n(y,\theta)|^2\ dy\ =\ 1\,.
\end{equation}

\vskip.2cm

\noindent
{\bf Assumption.}  
Fix $\theta_0\in [0,1[^N$, $n\in\NN$ 
and assume that $\lambda_n(\theta_0)$ is a {\bf simple} eigenvalue. 

\vskip.2cm

Define the associated nonnegative  frequency $\omega_n(\theta_0)$ satisfying the dispersion relation,
\begin{equation}
\label{disprel}
4\,\pi^2\omega_n^2(\theta)
\ =\ 
\lambda_n(\theta)\,.
\end{equation}
The group velocity is then defined by
\begin{equation}
\label{eq:groupvel}
\V \ :=\
-{\nabla_\theta\,\omega_n(\theta_0)}
\ = \ \frac{-\nabla_\theta\lambda_n(\theta_0)}{4\pi\sqrt{\lambda_n(\theta_0)}} \,.
\end{equation}
For any fixed $y\in\T^N$ the local speed of propagation
is given by 
$$
c(y)\  =\  \max_{1\leq j \leq N} \sqrt{\lambda_j(y)}
$$
where $\lambda_j(y)$ are the roots of the characteristic polynomial 
$$
p(y,\lambda) \ :=\  \det \left( A_0(y) - \lambda \rho_0(y) I \right) . 
$$
The maximal speed of propagation is 
$$
c_{max} \ :=\  \max_{y\in\T^N} c(y)\, .
$$

\begin{theorem}
\label{thm.speed}
The group velocity defined by \eqref{eq:groupvel} satisfies
$$
|\V| \ \leq\  c_{max} \,. 
$$
\end{theorem}

\begin{proof}
Introduce the operator 
$$
\Ao(\theta)\psi
\ :=\ 
- (\div_y +2i\pi\theta) \Big( A_0(y) (\grad_y +2i\pi\theta)\psi \Big) 
- \lambda_n(\theta) \rho_0(y) \psi \,.
$$
At the point $\theta_0$,
differentiate \eqref{celleq} with respect to $\theta$ in the direction 
of the covector $\xi$ to find
$$
\Ao(\theta) \xi.\nabla_\theta\psi_n = 
2i\pi \xi \cdot  A_0(y) (\nabla_y +2i\pi\theta)\psi_n
+ 2i\pi (\div_y +2i\pi\theta) (A_0(y) \xi \psi_n) + 
\xi.\nabla_\theta\lambda_n \rho_0(y) \psi_n .
$$
Taking the hermitian product of this inequality with $\psi_n$ 
yields
$$
\xi.\nabla_\theta\lambda_n (\theta) 
\ =\  2i\pi \int_{\T^N} \left( 
\psi_n A_0(y) \xi \cdot \overline{(\nabla_y +2i\pi\theta)\psi_n} 
- \overline{\psi_n} \xi \cdot  A_0(y) (\nabla_y +2i\pi\theta)\psi_n \right) dy \,.
$$
This implies the upper bound
$$
\left| \xi.\nabla_\theta\lambda_n (\theta) \right|
 \ \leq\ 
  4\pi 
\| \rho_0^{-1/2} A_0^{1/2} \|_{L^\infty(\T^N)}
\| \rho_0^{1/2} \psi_n \|_{L^2(\T^N)} 
\| A_0^{1/2} (\nabla_y +2i\pi\theta)\psi_n \|_{L^2(\T^N)} \,,
$$
which becomes
$$
\left| \xi.\nabla_\theta\lambda_n (\theta) \right|
\  \leq \
4\pi \,
c_{max}\, \sqrt{\lambda_n(\theta)}
$$
as desired.
\end{proof}

\begin{remark}   \rm The proof as given also proves a
bound of the propagation speed as a function of direction
by the corresponding fastest speeds of the original system
(see \cite{apr3} for an analogous result).
\end{remark}

\end{document}